\documentclass[11pt, oneside,a4]{article}   	
\usepackage{geometry}                		
\usepackage{graphicx}				

\usepackage{subeqnarray}

\title{An Novel Explicit Method to Solve Linear Dispersive Media for Finite Difference Time Domain Scheme}
\author{Hiroshi ABE (Three Wells) \footnote{habe@3wells-computing.com}}
\date{}							

\begin{document}
\maketitle
\abstract{A novel explicit method to model Lorentz linear dispersive media with finite difference method are presented.
The method shows an explicit method without any modification to the Leap-Frogging scheme.
The polarizations of the Lorentz media are convoluted simply in the current density. The numerical results are shown to be competitive accuracy with a conventional method. The codes used in the numerical investigation are available on the GitHub.}
\section{Introduction}
Solving dispersive media reactions to the electromagnetic fields are important in the plasma behaviour in metals.
The interaction is mutually affected so the phenomenon must be solve consistently.
Many algorithms for the problems have been introduced\cite{kelly1996},\cite{okoniewski2006} but many of them 
are very complicated to program.
A novel simple method that is very easy to implement is presented in the present paper.
\section{Equations}
\subsection{Maxwell's Equations}
The electromagnetic phenomena are governed by the Maxwell's equations,
\begin{eqnarray}
\frac{\partial \vec{B}}{\partial t} + \nabla \times \vec{E} & = & 0, \\
\frac{\partial \varepsilon \vec{E}}{\partial t} + \sigma \vec{E} + \vec{J} & = & \nabla \times \vec{B}/\mu, \label{eq:Amp}
\end{eqnarray}
where $\vec{E}, \vec{B}$ and $\vec{J}$ are the electric field, magnetic induction and current density, respectively. 
$\sigma$, $\mu$ and $\varepsilon$ are the conductivity, permeability and permittivity, respectively.
The macroscopic physical parameters, such as $\mu$, $\varepsilon$ and $\vec{J}$ are often governed by the other physics.
Concerning with dispersive media, several models are introduced and we have to solve the Maxwell's equations with the models.
Lorentz model is the one of the models.
\subsection{Lorentz Media}
The relative permittivity of the Lorentz model are described as the following equation.
\begin{equation}
\varepsilon(\omega) = \varepsilon_\infty 
+ \sum_{p=1}^P\frac{\Delta \varepsilon_p \omega_p^2}{\omega_p^2+2i\omega\delta_p-\omega^2}.
\end{equation}
We define,
\begin{equation}
\vec{P}_p(\omega) = \frac{\varepsilon_0 \Delta \varepsilon_p \omega_p^2}{\omega_p^2+2i\omega\delta_p-\omega^2} \vec{E}(\omega).
\label{eq:def_P}
\end{equation}
Inverse Fourier transforming Eq.(\ref{eq:def_P}) to obtain the governing equation in time domain,
\begin{equation}
\omega_p^2\vec{P}_p + 2\delta_p \frac{\partial \vec{P}_p}{\partial t} + \frac{\partial^2\vec{P}_p}{\partial t^2}
= \varepsilon_0 \Delta \varepsilon_p \omega_p^2 \vec{E}.
\label{eq:gov_P}
\end{equation}
With Eq.(\ref{eq:Amp}) to (\ref{eq:def_P}) we have modified Ampere-Maxwell equation,
\begin{equation}
\nabla \times \vec{B}/\mu_0 = \varepsilon_0 \varepsilon_\infty \frac{\partial \vec{E}}{\partial t} + \sigma \vec{E} + \vec{J}
+ \sum_{p=1}^P \frac{\partial \vec{P}_p}{\partial t}.
\label{eq:gov_Amp}
\end{equation}
We numerically solve Eq.(\ref{eq:gov_Amp}) and Eq.(\ref{eq:gov_P}) to evaluate the Lorentz media response exactly.
A novel method is introduced in the following section.
\section{Formulation with Transient Green Method}
Transient Green Method is a explicit but unconditionally stable method to solve transient wave equation\cite{habe1992}, and other type of equations \cite{abe2010stable}.
The method is to be shown applicable to the present problem.

We define an Impulse function which is defined as the following equation,
\begin{equation}
I(t, t_n;\Delta t) = \left\{ \begin{array}{ll} 
0 & \displaystyle t < t_n-\frac{\Delta t}{2}, \\
1 & \displaystyle t_n-\frac{\Delta t}{2} \leq t \leq t_n+\frac{\Delta t}{2}, \\
0 & \displaystyle t_n+\frac{\Delta t}{2} < t,
\end{array} \right.
\end{equation}
Let $G_p(t)$ be the Green function to the Eq.(\ref{eq:gov_P}) and we have,
\begin{eqnarray}
\left( \omega_p^2 + 2\delta_p \frac{\partial}{\partial t} + \frac{\partial^2}{\partial t^2} \right) G_p(t,t_n;\Delta t)
& = & I(t,t_n;\Delta t). \label{eq:lorentz_g} \\
& & t\ge t_n+\Delta t/2, \nonumber
\end{eqnarray}
The equation is linear equation so we can obtain the polarization vector $\vec{P}_p$ by,
\begin{equation}
\vec{P}_p(t) = \varepsilon_0 \Delta \varepsilon_p \omega_p^2 \sum_{n=0}^N \vec{E}^n G_p(t,t_n;\Delta t),
\label{eq:green_sum}
\end{equation}
Eq.(\ref{eq:green_sum}) gives polarization vector $P$ at time $t$. 
But the formulation is quite inconvenient because you have to sum all terms whenever the time step is advanced.
To avoid this inconvenience a trick is introduced to obtain an explicit recurrence formulation to evaluate the current time values from values of one time step before only.

First of all, the Green function $G_p$ with respect to the Impulse function $I$ is to be solved.
Impulse function $I(t,t_n;\Delta t)$ is Fourier transformed to $I(\omega;t_n, \Delta t)$,
\begin{eqnarray}
I(\omega;t_n,\Delta t) & = & \int_{-\infty}^\infty I(t,t_n;\Delta t) e^{-i\omega t} dt, \nonumber \\
& = & \int_{t_n-\frac{\Delta t}{2}}^{t_n+\frac{\Delta t}{2}} e^{-i\omega t} dt, \nonumber \\
& = & \frac{e^{i\omega \Delta t/2} - e^{-i\omega \Delta t/2}}{-i\omega} e^{-i\omega t_n} \label{eq:rect_f}
\end{eqnarray}
Substituting Eq.(\ref{eq:rect_f}) to Eq.(\ref{eq:lorentz_g}) and the Green function is obtained as,
\begin{eqnarray}
G_p(\omega;t_n,\Delta t) & = & \frac{I(\omega;t_n,\Delta t)}{\omega_p^2 +i\omega 2\delta_p - \omega^2}, \nonumber \\
& = & \frac{e^{i\omega \Delta t/2} - e^{-i\omega \Delta t/2}}
{i\omega( \omega^2-i\omega 2\delta_p -\omega_p^2)} e^{-i\omega t_n}.
\end{eqnarray}
The Green function $G_p(t;t_n,\Delta t)$ can be obtained by Fourier inverse-transformation to the $G_p(\omega;t_n,\Delta t)$.
\begin{equation}
G_p(t;t_n,\Delta t) = \frac{1}{2\pi} \int_{-\infty}^\infty G_p(\omega; t_n,\Delta t) e^{i\omega t} d\omega.
\end{equation}
The integration can be solved by the Cauchy's integral theorem.
There are three poles of order 1 at $\omega = 0, z^{\pm}$, where
$z^\pm = i\delta_p \pm \sqrt{\omega_p^2-\delta_p^2}$.
\begin{eqnarray}
G_p(t;t_n,\Delta t) & = & \frac{2\pi i}{2\pi} \sum Res(0, z^\pm), \nonumber \\
& = & i \left( \frac{e^{iz^+ \frac{\Delta t}{2}} - e^{-iz^+ \frac{\Delta t}{2}}}
{iz^+(z^--z^+)} e^{iz^+ (t-t_n)}
+ \frac{e^{iz^- \frac{\Delta t}{2}} - e^{-iz^- \frac{\Delta t}{2}}}
{iz^-(z^+-z^-)} e^{iz^- (t-t_n)}
\right). \label{eq:sol_G}
\end{eqnarray}
Substituting Eq.(\ref{eq:sol_G}) to Eq.(\ref{eq:green_sum}) to obtain,
\begin{eqnarray}
\vec{P}_p(t) & = & \varepsilon_0\Delta\varepsilon_p \omega_p^2 \sum_{n}  \vec{E}^n G_p(t;t_n,\Delta t), \\
& = & \varepsilon_0\Delta\varepsilon_p \omega_p^2  \sum_{n=0}^N \vec{E}^n\nonumber \\
& \times & \left( \frac{e^{iz^+ \frac{\Delta t}{2}} - e^{-iz^+ \frac{\Delta t}{2}}}
{z^+(z^--z^+)} e^{iz^+ (t-t_n)}
+ \frac{e^{iz^- \frac{\Delta t}{2}} - e^{-iz^- \frac{\Delta t}{2}}} 
{z^-(z^+-z^-)} e^{iz^- (t-t_n)} \right) \nonumber \\
& = & \varepsilon_0\Delta\varepsilon_p \omega_p^2\left[ e^{iz^+(t-t_N)} \vec{F}_N^+ + e^{iz^-(t-t_N)} \vec{F}_N^-
\right],
\end{eqnarray}
where,
\begin{eqnarray}
\vec{F}_N^\pm = \vec{F}_{N-1}^\pm \cdot e^{iz^\pm(t_N-t_{N-1})} + 
 \frac{e^{iz^\pm \frac{\Delta t}{2}} - e^{-iz^\pm \frac{\Delta t}{2}}}{z^\pm(z^\mp-z^\pm)} \vec{E}^N \label{eq:recur_F}.
\end{eqnarray}
Fuerthermore, we can obtain $\partial \vec{P}_p(t)/\partial t$ by simple multiplications.
\begin{equation}
\frac{\partial \vec{P}_p(t)}{\partial t} 
= i\varepsilon_0\Delta\varepsilon_p \omega_p^2\left[ z^+ e^{iz^+(t-t_N)} \vec{F}_N^+ + z^- e^{iz^-(t-t_N)} \vec{F}_N^-
\right]
\label{eq:sol_Pt}
\end{equation}
Eq.(\ref{eq:recur_F}) gives an explicit scheme to solve Eq.(\ref{eq:gov_P}).
Furthermore we can solve Eq.(\ref{eq:gov_Amp}) without any modification to the Leap-Frog scheme for solving the Maxwell's equations using Eq.(\ref{eq:sol_Pt}) with the same as the current density term.

\section{Numerical Investigation}
\subsection{Model}
\begin{figure}[htbp]
\begin{center}
\includegraphics[width=12cm]{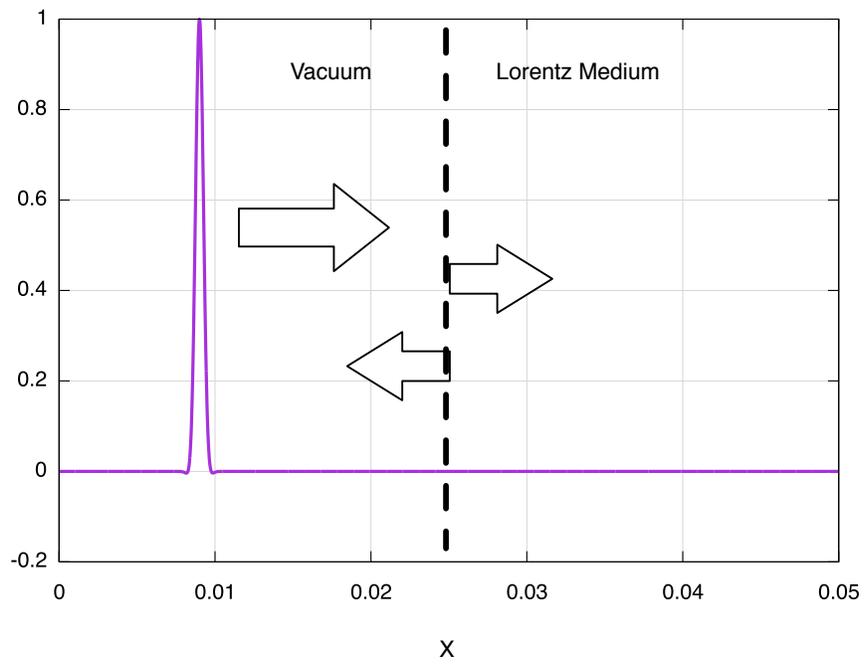}
\caption{{\bf Numerical model schematic view.}}
The boundary of the mediums, vacuum area and Lorentz medium area, is located at the centre of the numerical space.
The Gaussian impulse is applied from the left boundary.
\label{fig:model}
\end{center}
\end{figure}

We suppose one dimensional case with vacuum space and a Lorentz dispersive half-space.
The boundary is located at the centre of the simulation space as is shown in the Fig.(\ref{fig:model}).

A Gaussian pulse,
\begin{equation}
E_{boundary}(t) = e^{-\frac{(t-t_0)^2}{2\Delta T^2}} \cos \omega_0 (t-t_0),
\end{equation}
is applied at the left boundary.
The boundaries are set to absorbing boundary condition of Mur's first-order scheme but the simulation space is set to be long enough to avoid unnecessary reflections at the boundaries.
We probe the electric field at three points to sample the time series of the filed values.
We calculate the reflection coefficients from the time series of simulations and compare them to the exact refection coefficients.
\subsection{Reflection Coefficient}
The reflection coefficient $R$ of a Lorentz dispersive medium with respect to a vacuum space is given by the following equation\cite{jackson}.
\begin{equation}
R = \frac{E_{ref}}{E_{inc}} = \frac{\sqrt{\frac{\mu_0\varepsilon_0 \varepsilon'}{\mu'\varepsilon_0}}-1}{\sqrt{\frac{\mu_0 \varepsilon_0 \varepsilon'}{\mu'\varepsilon_0}}+1}.
\end{equation}
For $\mu_0=\mu'$ we have,
\begin{equation}
R(\omega) = \frac{\sqrt{\varepsilon'(\omega)}-1}
{\sqrt{\varepsilon'(\omega)}+1}.
\end{equation}

\subsection{Numerical Simulation}
The one dimensional model is simulated with TGM here.
Also Auxiliary Differential Equation Method (ADEM) is performed as comparison.
The two numerical results are compared to the exact solution.

The one dimensional (x-axis) Maxwell's equations to solve are,
\begin{eqnarray}
\frac{\partial B_z}{\partial t} & = & -\frac{\partial E_y}{\partial x}, \\
\frac{\partial (\varepsilon_\infty E_y + P_y)}{\partial t} & = & \frac{\partial B_z}{\partial x}.
\end{eqnarray}
We need $\frac{\partial \vec{P}_p(t)}{\partial t} $ at $t=t_N+\frac{\Delta t}{2}$ so,
\begin{equation}
\frac{\partial P_{yp}(t_N+\frac{\Delta t}{2})}{\partial t} = 
i\varepsilon_0\Delta\varepsilon_p \omega_p^2\left[ z^+ e^{iz^+\frac{\Delta t}{2}} F_{ypN}^+ + z^- e^{iz^-\frac{\Delta t}{2}} F_{ypN}^-
\right].
\end{equation}

The numerical parameters are shown in Table \ref{tab:params}.
\begin{table}[ht]
\caption{\bf Parameters for Numerical Simulation}
\begin{center}
\begin{tabular}{|c|c|} \hline
Parameters & Values (MKSA) \\ \hline \hline
System Length (L) & 0.05 \\ \hline
Grids (N) & 3000 \\ \hline
Light Speed (C) & $2.99792458\times10^8$ \\ \hline
$\mu_0$ & $4\pi \times 10^{-7}$ \\ \hline
$\varepsilon_0$ & $1.0 /(\mu_0C^2)$ \\ \hline
$\Delta x$ & L/(N-1) \\ \hline
$\Delta t$ & $0.9\Delta x/C$ \\ \hline
$\omega_0$ & $2\pi 100\times 10^9$ \\ \hline
$\Delta T$ & $1.0 \times 10^{-12}$ \\ \hline
$t_0$ & $1.0 \times 10^{-11}$ \\ \hline
$\varepsilon_\infty$ & 1.5 \\ \hline
$\Delta \varepsilon$ & 3.0 \\ \hline
$\omega_p$ & $2\pi\times20.0 \times 10^9$ \\ \hline
$\delta_p$ & $0.1\omega_p$ \\ \hline
$\sigma$ & 0.0 \\ \hline
\end{tabular}
\end{center}
\label{tab:params}
\end{table}%

\subsection{Numerical Result}
\begin{figure}[htbp]
\begin{center}
\input{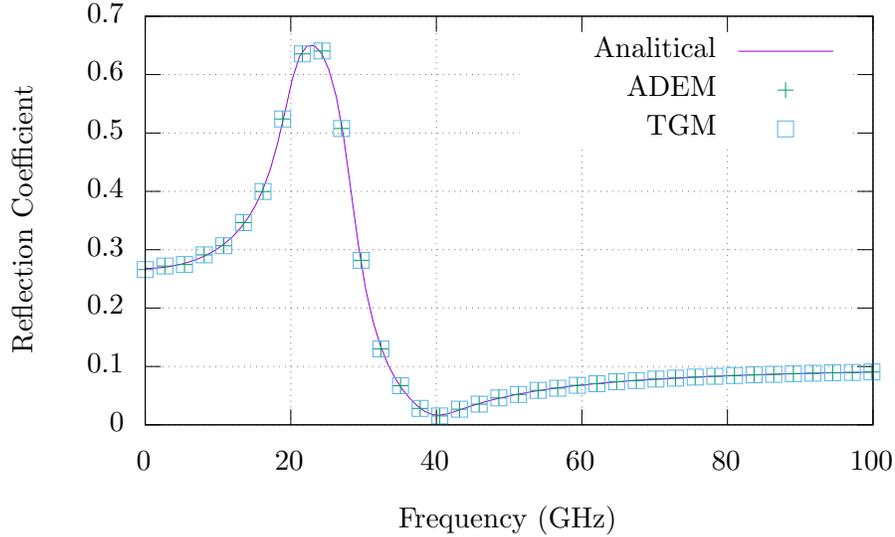}
\label{fig:reflection}
\end{center}
\caption{\bf Refection Coefficients}
All the resultant coefficients are well matched. The TGM shows the equivalent performance to the ADEM.
\end{figure}

The reflection coefficients in the frequency domain are shown in Fig.\ref{fig:reflection}.
As a comparison, the conventional ADE method is also shown in the figure.

\section{Summary}
A simple method for evaluating the linear dispersive media is presented in the paper.
It is very simple to program with normal FDTD Leap Frog method without any modification to the Leap Frogging.
It has second ordered temporal accuracy and give competitive results to the conventional ADEM.
The proposed method should apply to the other linear dispersive media, such as Debye, Drude and the others.

\appendix
\section{Programs}
All the programs are available from the GitHub \cite{github}.
They are Python scripts and some are Gnuplot scripts.

\end{document}